# Statistical thinking:
# From Tukey to Vardi and beyond

Larry Shepp[1],*

*Rutgers University*

**Abstract:** Data miners (minors?) and neural networkers tend to eschew modelling, misled perhaps by misinterpretation of strongly expressed views of John Tukey. I discuss Vardi's views of these issues as well as other aspects of Vardi's work in emision tomography and in sampling bias.

Statistics is not in my main skill set but I take this opportunity to record a few thoughts about it I have accumulated over the years.

In the '60's, John Tukey and his followers brought exploratory data analysis into statistics, partly as a revolt against what was then perceived as an overly rigid and brittle mathematical modelling philosophy that held sway at that time. Some problems seemed to demand such a purely data-driven approach where data mining methods in the absence of mathematical modelling is the driving philosophical methodology. One did not want to be biased by preconceived ideas about the origin of the data by formulating a model but instead allowed the data to "speak for itself". Vardi liked mathematical modelling and was very good at it. He also promoted data mining, depending on the problem and thus straddled both philosophies. He and I often debated these issues, and were often in friendly disagreement.

I will try to argue with concrete examples of work of Vardi and others in statistics that the pendulum should again swing back a bit towards encouraging more mathematical modelling to obtain maximal benefit from the use of statistical procedures by allowing physics, biology, and other fields of science to enter the statistical problem formulation via mathematical modelling of the specific statistical problem at hand. I would argue that the solution to a specific problem ought to somehow depend on the problem itself, which is not the case with neural-nets and other data-driven approaches that live mostly or entirely within the data or training set of the problem.

Data-driven statistics has the danger of isolating statistics from the rest of the scientific and mathematical communities by not allowing valuable cross-pollination of ideas from other fields. To illustrate these ideas I will discuss among other concrete examples of statistics problems: emission tomography, machine learning, sampling bias. These topics were debated frequently by Vardi and me. I will do my best to give Vardi's side as honestly as possible. Needless to say, I wish he were here to continue the debate. I will quote Tukey and/or Vardi on issues I will raise, and you should be aware that people who quote absentees don't allow the quotees to modify the positions they are being quoted on unless they are in agreement with the

---

*Research partially supported by National Science Foundation Grant DMS-0504387.

[1]Department of Statistics and Biostatistics, Rutgers University, Hill Center, Busch Campus, Piscataway NJ 08854, USA, e-mail: shepp@stat.rutgers.edu

*Keywords and phrases:* Data mining, neural nets, statistical modelling, emission tomography, sampling bias.

*AMS 2000 subject classification:* 6207.

quoter's own positions. However, if I quote the current views of Tukey and Vardi inaccurately it is inadvertent!

## 1. Emission tomography and the EM algorithm

Emission tomography has the advantage over CT scanning that it can be used to study metabolism. Thus one can in principle learn, say, where in the brain, higher cognition is taking place by tagging psychoactive substances with radioactivity. These substances move to the part of the brain which is active during the performance of a higher cognition task, and then radiate $\gamma$ rays which are then detected. Vardi and I wrote [1] on ET and then he wrote [2] and several others without my participation:

His role was large indeed: I was trying to find a maximum likelihood estimate of the unknown emission density because one could write down in closed form the likelihood of the given detected counts for an arbitrary emission density using the well-known Poisson model for radioactive emissions from a given emission density. Since emission CT is "photon-limited" it is reasonable to use maximum likelihood estimation as a driving modelling approach. Vardi was familiar with the EM algorithm, with the notion of "over-fitting", and also with Kullback-Leibler theory and alternating convergence techniques. I was indeed lucky that Vardi was at Bell Labs at the time and that he got very interested in this project. Our second paper above, which is much better written than the first, was entirely his and Linda Kaufman's doing.

Later on, with the advent of functional MRI, I felt that brain physiology is better modelled by the physics of MRI. The reason is that emission CT is too slow to study rapid events in the brain which take place in a fraction of a second. Statistics still plays a major role in fMRI brain physiology but the methodology of the EM algorithm is no longer involved because the physics is now different and there is no likelihood to maximize in fMRI. I feel that one should employ methods that reflect the physics of the problem at hand rather than the methods one happens to know. I said this so often in our discussions that Vardi used to laugh at me.

Vardi turned his attention to problems in which he could use parallel EM type methodology to solve reconstruction problems of similar type but without having the clear basis in physics for setting a criterion for a maximization. We generally agreed to disagree on this methodological point. I must admit that he got some very nice results on various problems just by using convergence methods of EM type even with no basis in physics. Vardi showed for example how to reconstruct the traffic on each leg of a graphical network from the overall traffic between the nodes of the network, which he called "network tomography" and which is indeed strongly analogous to emission CT. This debate between Vardi and me became even more relevant in the next topic. Of course it was always a friendly debate though sometimes I got too loud.

## 2. What would the founding fathers have said about neural nets and machine learning?

I often tease the neural net community by asking them to design a neural net that would take the CT, ET, or MRI reconstruction and try to find a tumor in it or decide that no tumor is present. This might be possible of course and indeed there are many papers written on this topic which use image enhancement methods to delineate blobs that a radiologist could then make a decision as to whether to suspect a tumor which would save effort, presumably. But to make things more interesting, I dared the neural net community to work not with the CT, ET, or MRI *reconstruction*, but instead with the *raw or measured data*. After all if neural net technology really works, who needs Radon inversion or Fourier transforms; why not use it on the raw data directly? I made this challenge to point out the mindlessness of throwing mathematical modelling away completely.

Another difficult problem area, which is more difficult to use mathematics on and in which there are many attempts to use neural nets or other statistical approaches to pattern recognition problems is on the problem of automated understanding of human speech. I must admit that much progress has been made on it by Allen Gorin who used his ideas about "salience" to build systems like "How May I Help You" for AT&T which involves direct discourse with a human being who makes one of a limited number of requests from an operator or a business office. Salience led quickly to a working system which accomplished this limited task. It avoids the need for the usual confusion of a system using things like "if you want ... press one, if ... press two" which irks many people no end. Here machine learning avoids mathematical modelling or true understanding and yet it does accomplish the job. Vardi often would rub my nose in this fact.

At the other end of the spectrum of machine understanding of speech would be a problem such as machine translation of Russian text into English. Here decision theory or salience can hardly be expected to find an appropriate $L^2$-loss function to be used to train an algorithm to converge to a good translator. It seems clear that when we do construct a robotic translator it will have to somehow understand the "meaning" of the passage to be translated and it will no doubt use some modular system based on key word and rules rather than on the standard statistical decision theory, salience, or neural net approaches.

What was Vardi's opinion of this controversy? Did he agree with Tukey's position, with mine, or with neither? I'm not sure; probably neither, but I have to say that he liked to reuse the methods of statistics with no basis in modelling more than I do, and did it occasionally more successfully than I expected. He used clever convergence methods and did not require there to be any physics behind the model. I would often argue that this gives no way for mathematics and physics to interact with statistics. I think he agreed with me at least once, but I am sure both Tukey and Vardi agree completely with my position today!

There seems to me to be a disturbing trend by statisticians to use "standard statistical methodology" to solve problems which may not be amenable to simple approaches. Let me give a few examples of such problems and let's try to guess what John Tukey and Yehuda Vardi's opinions would be on the issues.

Robotics is an important problem to solve for the future development of society. Can statistics play a key role in this important area? A useful robot would relieve us of the need to perform routine tasks and would also provide entertainment and companionship when desired. Most problems in robotics call for the development of algorithms for automated pattern recognition often called machine learning as we discussed above. Examples of robotic tasks include recognizing hand-written characters so that addresses on envelopes can be automatically sorted. More profound pattern rocognition tasks include speech understanding so that a robot could respond to us directly and converse with us in order to take instructions and perhaps provide companionship. Chess playing robots are already well-developed and provide an opponent at any time for people to play with.

Another conceptual problem I enjoy contemplating is to program a computer to recognize whether a given picture of an animal is that of a dog or that of a cat. Of course small children can do this with high accuracy in most cases but it may not be so easy to write such a program without some understanding of what is the "salient feature" that is the real difference between the two animals. Is it reasonable to try to find such a salient feature by a neural net on a training set or is this likely to find some "feature" that has little or nothing to do with how a child actually does it? Instead the neural net might just find some commonality of all the dogs in the training set that is really totally irrelevant to the problem. Also it seems likely that the system of the neural net would not be modular since there would not be any way to build on it to include recognition of other animals but instead one would have to start anew.

While I would argue that each particular problem of pattern recognition should bring a different solution based on some feature from the real world that we can see is relevant and then find an algorithm to look for the feature, I think Tukey might argue just the opposite. After all, he advocated looking at the data of a particular statistical problem without trying to model it, so as not to prejudice yourself with preconceived ideas. Exploratory data analysis is similar to the use of salient features found automatically via a neural net on a training set except that Tukey himself or some other data miner would be doing the analysis rather than a program. But Tukey's articulate urging made neural nets more attractive and so in a sense Tukey may have stimulated neural nets. Would he have changed his mind today in the light of some of the overblown claims of neural nets? I think so, or at least I hope so.

**Character recognition; an easy machine learning problem.** The problem of robotics has statistical components because machine learning is based on training set data which has lots of randomness. Thus a neural net for the "post office problem" uses a training set of tens of thousands of examples of hand-written characters where we are given which particular character was written by a person as well as a (say $16 \times 16$ zero-one) digitized image of the character in each example. The algorithm classifies each image and names a character. One statistical approach is to minimize the $L^2$ distance between the new image and the set of all images for a given character. The neural net approach is similar: they use a scheme where a collection of linear functionals is used and training takes place to choose automatically which linear functional is most salient on the training set. This approach has the advantage that one does not have to think, and one does not have to use any physics or any preconceived notions to guide the convergence of the neural net to the *meaningless* salient feature. In this case, however, the problem is not all that hard and almost any scheme can recognize characters with an error rate of a few percent.

While Tukey, Vardi, and I were at Bell Labs there was much work being done on the above problem of recognizing just handwritten characters from 0 to 9 for automated zipcode reading for the post office sorting application. Let me relate my recollections as an outside observer of this effort which involved several large groups of people. One of the engineers, Patrice Simaud, suddenly announced a great advance and Trevor Hastie, who had been working on the problem with statistical methods, and had been getting less successful error rates, courageously invited Patrice to present his results in our statistics seminar. Instead of error rates in the region of 2 or 3 percent obtained by the methods of the neural nets and the statistical loss-function approach, Patrice announced an algorithm that was virtually error-free except in those cases where even a human being couldn't recognize a character that was sloppily written. How did Patrice accomplish this?

He used "physics" a modelling of the problem in a clever way. Of course this problem is not as hard as other automated learning problems mentioned above but it still is hard enough to be done incorrectly if one does not think and decides to use statistics or neural nets blindly.

Patrice took advantage of the fact that there are mathematical invariances present in the problem. The image of the number being written will depend on how much ink is in the pen, or how thick the pencil point is, how the envelope is oriented relative to the wrist, etc. This type of inputting something from the real world to solve the problem is what I like to see. Patrice introduced 7 transformations of the image of the character to be identified. Each character was thickened a bit (as if it were written by a pen with a larger tip) to get transformation one, was rotated to get transformation two, was translated in two directions to get transformations three and four, and so on. I don't recall the other 3 transformations he used. This gave him then 7 new points in $16 \times 16 = 256$ dimensional space and determined together with the original image a 7 dimensional hyperplane in 256 dimensional space. He then dropped a perpendicular onto this plane from each of the 10,000 training characters and identified the new character as that of the closest point to the hyperplane in the training set. He had a method that completely solved the problem. Clever!

The sad part of the story was that he presented his solution to our seminar as a neural net method which it was certainly not. Patrice thought about the problem from a modelling point of view and looked for a way to let the solution depend on the problem. Since Patrice's approach yielded a perfect algorithm one would think that it would be the preferred method today but it is not! The method had a drawback: it is slow. Dropping the 10,000 perpendiculars took a lot of time. One might think that in view of the great improvement in the error rate, one might put effort into computational speed-up of his method, but instead everyone worked still harder on the neural net approach. After considerable effort, the neural net community found that a fast algorithm could be obtained to find salient features which matched the performance of Patrice's method. I would venture to guess that there was some guiding of the convergence of the the neural net to incorporate some features that were understood only after seeing how Patrice's algorithm performed. Alas, even today even Patrice himself insists that neural nets is the better way to go for character identification. Argh!

## 3. Sampling bias

Vardi was a great statistician and would instinctively move to a simple and clear understanding of how to think about a new statistical problem. Often his way of thinking would involve little mathematical technique and so others could copy his ideas easily and since he was always very open and generous, he sometimes lost some credit for his insightful ideas. One exception was the area of sampling bias where he was <u>the</u> recognized expert. Others at this conference will characterize Vardi's work on sampling bias better than I can, but he and I did solve one cute mathematics problem where the motivation originated in sampling bias. We wrote [3] with Ben Logan.

The problem posed by Vardi was to find the class of all cdF's, $F = F(x)$, for the lifetime of a lightbulb for which the residual lifetime would be a scaled version of $F$, $F(qx)$. Vardi wanted to know whether there was an analog to the result for the exponential, $F(x) = 1 - e^{-x}$, that the residual lifetime cdF is the same as the lifetime cdF itself, for other values of $q$. We showed that the class, $C_q$, of all $F$'s with Vardi's scaling property is either empty (if $q > 1$) unique (if $q = 1$), or uncountable and convex if $(0 < q < 1)$. However even for $q = .5$ there is essentially only a single cdF in $C_q$ in the sense that all the graphs of $F \in C_{.5}$ would fit inside a single pencil-line curve, as we showed. The class $C_q$ gets large as $q \to 0$.

Sampling bias plays its role in lifetime statistics as Vardi emphasized by observing that if one samples lightbulbs or obituaries in the NYTimes one gets an incorrect sampling of performance of individuals since one is looking mostly at *long-lived* ones.

Vardi had great taste for what constituted a good mathematics problem. We all will miss his insights, his leadership, and his great sense of humor.